\documentclass[11pt]{amsart}

\usepackage{amsmath}
\usepackage{amssymb}
\usepackage{graphicx}
\usepackage[hyphens]{url} 
\usepackage[utf8]{inputenc}
\usepackage{mathtools}  
\usepackage{amsthm}
\usepackage{color}
\usepackage{nicefrac}
\usepackage{tcolorbox}

\newtheorem{theorem}{Theorem}[section]
\newtheorem{corollary}[theorem]{Corollary}

\theoremstyle{definition}
\newtheorem{definition}[theorem]{Definition}
\newtheorem{example}[theorem]{Example}
\newtheorem{remark}[theorem]{Remark}

\numberwithin{equation}{section}

\title[On local uniqueness of normalized Nash equilibria]{On local uniqueness \\ of normalized Nash equilibria}

\author[V. Shikhman]{Vladimir Shikhman}
\address[V. Shikhman]{Department of Mathematics, Chemnitz University of Technology, Reichenhainer Str. 41, 09126
Chemnitz, Germany}
\email{\tt vladimir.shikhman@mathematik.tu-chemnitz.de}



\begin{document}

\begin{abstract}
For generalized Nash equilibrium problems (GNEP) with shared constraints we focus on the notion of normalized Nash equilibrium in the nonconvex setting. The property of nondegeneracy for normalized Nash equilibria is introduced. Nondegeneracy refers to GNEP-tailored versions of linear independence constraint qualification, strict complementarity and second-order regularity. Surprisingly enough, nondegeneracy of normalized Nash equilibrium does not prevent from degeneracies at the individual players' level. We show that generically all normalized Nash equilibria are nondegenerate. Moreover, nondegeneracy turns out to be a sufficient condition for the local uniqueness of normalized Nash equilibria. We emphasize that even in the convex setting the proposed notion of nondegeneracy differs from the sufficient condition for (global) uniqueness of normalized Nash equilibria, which is known from the literature.        

\end{abstract}

\maketitle


\section{Introduction}

We consider the class of generalized Nash equilibrium problems (GNEP). The term "generalized" refers to the fact that not only the objective functions of the players, but also their feasible sets mutually depend on other players' strategies. This dependence is explicitly given by means of the so-called shared constraints. They model the common limitation of resources, e.g. in form of the budget constraint or the common use of transportation facilities, see \cite{faccinei:2010}. The class of GNEPs with shared constraints has been widely studied in the literature so far. The starting point is the seminal paper \cite{rosen:1965}, where Nash equilibria for GNEPs were examined. It turns out that the presence of shared constraints makes the set of Nash equilibria look geometrically complex. In general, we cannot expect that it consists of isolated, not to speak of unique points. Rather than that, the set of Nash equilibria exhibits higher-dimensional parts, kinks and boundary points, see \cite{dorsch:2013}. In order to avoid this non-uniqueness drawback, the concept of normalized Nash equilibrium has been introduced in \cite{rosen:1965}. Normalized Nash equilibria are characterized by fixing the ratios of multipliers for the shared constraints in the corresponding system of concatenated Karush-Kuhn-Tucker conditions for the players' optimization problems. In the convex setting, fundamental results on existence and uniqueness of normalized Nash equilibria were provided in \cite{rosen:1965}. The study of normalized Nash equilibria has been continued in \cite{faccinei:2007} and \cite{nabetani:2011} by characterizing them as solutions of  parametrized variational inequalities. Another possibility to compute normalized Nash equilibria was suggested in \cite{heusinger:2012}, where a fixed-point formulation was used for this purpose. An interesting generalization of the notion of normalized Nash equilibrium, called restricted Nash equilibrium, was presented in \cite{fukushima:2011}.

In this paper, our aim is to study the local uniqueness of normalized Nash equilibria in the nonconvex setting. For that, we introduce the concept of nondegeneracy for normalized Nash equilibria. Nondegeneracy refers to GNEP-tailored versions of linear independence constraint qualification, strict complementarity and second-order regularity, see Definition \ref{def:nd-m-st}. We show that GNEP-LICQ holds at all feasible points of a generic GNEP, see Theorem \ref{thm:licq-gen}. Moreover, all normalized Nash equilibria are shown to be nondegenerate in a generic sense, i.e. for a dense subset of GNEP defining functions with respect to a suitable topology, see Theorem \ref{thm:genericI}. We emphasize that nondegeneracy of normalized Nash equilibrium does not prevent from degeneracies at the individual players' level. It may well happen that at a nondegenerate normalized Nash equilibrium e.g. the standard LICQ or the standard strict complementarity fail to hold for players' parametric optimization problems, see Example \ref{ex:deg}. Our main result says that nondegenerate normalized Nash equilibria are nevertheless locally unique, see Corollary \ref{cor:unique}. In other words, nondegeneracy is a sufficient condition for the local uniqueness of normalized Nash equilibria. We point out that the proposed notion of nondegeneracy differs from the sufficient condition for (global) uniqueness given in \cite{rosen:1965}. Surprisingly enough, this is the case even in the convex setting, see Example \ref{ex:compar}.        

The paper is organized as follows. Section \ref{sec:nne} is devoted to the notion of normalized Nash equilibrium and preliminaries. In Section \ref{sec:lu}, we address the local uniqueness of normalized Nash equilibria by introducing the concept of their nondegeneracy.

Our notation is standard. The cardinality of a finite set $A$ is denoted by $|A|$. The $n$-dimensional
Euclidean space is denoted by $\mathbb{R}^n$ with the coordinate vectors $e_i,i= 1,\ldots, n$. Given a twice continuously differentiable function $f:\mathbb{R}^n\rightarrow \mathbb{R}$, $D f$ denotes its gradient as the row vector, and $D^2f$ stands for its Hessian. Given a continuously differentiable function $F:\mathbb{R}^n\rightarrow \mathbb{R}^m$, $D F$ denotes its Jacobian $(m,n)$-matrix. 


\section{Normalized Nash equilibrium}
\label{sec:nne}
The generalized Nash equilibrium problem (GNEP) can be
stated as follows. Each player $\nu$ from a finite set $\mathcal{N} =\{1,\ldots, N\}$ solves the following parametric nonlinear programming
problem:
\[
  P_\nu(x^{-\nu}): \quad \min_{x^\nu \in {X}^{\nu}} f^\nu(x^\nu, x^{-\nu}) \quad \mbox{s.t.} \quad G_j(x^\nu, x^{-\nu}) \geq 0, \quad j \in \mathcal{J},
\]
where $\mathcal{J}=\{1,\ldots,|\mathcal{J}|\}$ is the index set of shared constraints, and the player $\nu$’s strategy set is given by
\[
   X^\nu = \left\{x^\nu \in \mathbb{R}^{n_\nu} \, \left| \, g_j^\nu(x^\nu) \geq 0, j \in \mathcal{J}^\nu\right.\right\},
\]
where $\mathcal{J}^\nu=\{1,\ldots,|\mathcal{J}^\nu|\}$ is the index set of individual constraints.
We assume that all defining functions $f^\nu, g_j^\nu, j \in \mathcal{J}^\nu, \nu \in \mathcal{N}$, $G_j, j \in \mathcal{J}$, are twice continuously differentiable.
By means of $x^{-\nu}$ we denote as usually the vector formed by all the players’ variables except those of player $\nu$. For $x \in \mathbb{R}^n$ with $n=\sum_{\nu=1}^{N} n_\nu$ we occasionally write $(x^{\nu}, x^{-\nu})$ in order to emphasize the $\nu$-th player’s variables within $x$, or simply $(x^1,\ldots, x^N)$.

For short, we set
\[
   f=\left(f^\nu, \nu \in \mathcal{N}\right),\quad
   g=\left(g^\nu_j, j \in \mathcal{J}^\nu, \nu \in \mathcal{N}\right),\quad
   G=\left(G_j, j \in \mathcal{J} \right),
\]
and associate the class of GNEPs with their defining functions: 
\[
   (f,g,G) \in C^2\left(\mathbb{R}^n,\mathbb{R}^N\right)\times \prod_{\nu=1}^{N} C^2\left(\mathbb{R}^{n_\nu},\mathbb{R}^{\left|\mathcal{J}^\nu\right|}\right) \times C^2\left(\mathbb{R}^n,\mathbb{R}^{\left|\mathcal{J}\right|}\right).
\]
The latter product space will be endowed with the strong (or Whitney) $C^2$-topology, denoted by $C^2_s$, see e.\,g. \cite{hirsch:1976}. The $C^2_s$-topology is generated by allowing perturbations of the functions, their gradients and Hessians, which are controlled by means of continuous positive functions. 
We say that a property is generic for GNEP if there exists a $C_s^2$-open and -dense subset $\mathcal{D} \subset C^2\left(\mathbb{R}^n,\mathbb{R}^N\right)\times \prod_{\nu=1}^{N} C^2\left(\mathbb{R}^{n_\nu},\mathbb{R}^{\left|\mathcal{J}^\nu\right|}\right) \times C^2\left(\mathbb{R}^n,\mathbb{R}^{\left|\mathcal{J}\right|}\right)$, such that all GNEPs with defining functions $(f,g,h) \in \mathcal{D}$ fulfill this property.

The feasible set of GNEP is given by
\[
   X= \left\{\left. x \in \prod_{\nu=1}^{N}X^\nu \, \right| \,G_j(x^1,\ldots, x^N) \geq 0, j \in \mathcal{J}\right\}.
\]
The index set of the active shared constraints at a GNEP feasible point $\bar x \in X$ is denoted by
\[
   \mathcal{J}_0(\bar x) = \left\{j\in \mathcal{J} \,\left|\, G_j(\bar x)=0\right.\right\}.
\]
The index set of the active player $\nu$'s individual constraints at $\bar x^\nu \in X^\nu$ is analogously denoted by
\[
   \mathcal{J}_0^\nu(\bar x^\nu) = \left\{j\in \mathcal{J}^\nu \,\left|\, g_j^\nu(\bar x^\nu)=0\right.\right\}.
\]

Let us start with the celebrated notion of Nash equilibrium.

\begin{definition}[Nash equilibrium]
\label{def:nash}
A GNEP feasible point $\bar x \in X$ is called Nash equilibrium if
$\bar x^\nu$ solves $P_\nu(\bar x^{-\nu})$ for each $\nu \in \mathcal{N}$.
\end{definition}


It has been observed already in \cite{rosen:1965} that even in the convex setting GNEPs usually have multiple Nash equilibria  due to the presence of shared constraints. This motivated Rosen to study their subclass, called normalized Nash equilibria.  

\begin{definition}[Normalized Nash equilibrium, \cite{rosen:1965}] 
\label{def:norm-nash}
Let us fix the positive parameters $r=\left(r^\nu, \nu \in \mathcal{N}\right)$.
A Nash equilibrium $\bar x \in X$ is called normalized if there exist real multipliers $\bar \lambda^\nu_j$, $j \in \mathcal{J}^\nu$, and $\bar \Lambda^\nu_j$, $j \in \mathcal{J}$, $\nu \in \mathcal{N}$, such that for every $\nu \in \mathcal{N}$ it holds:
\begin{equation}
    \label{eq:kkt1}
    D_{x^\nu} f^\nu(\bar x) = \sum_{j \in \mathcal{J}^\nu} \bar \lambda^\nu_j D g_j^\nu(\bar x^\nu)+\sum_{j \in \mathcal{J}} \bar \Lambda^\nu_j D_{x^\nu} G_j(\bar x),
\end{equation}
\begin{equation}
    \label{eq:kkt2a}
    \bar \lambda^{\nu}_j g_j^\nu(\bar x^\nu) =0, \bar \lambda^{\nu}_j \geq 0, g_j^\nu(\bar x^\nu) \geq 0, \quad j \in \mathcal{J}^\nu,
\end{equation}
\begin{equation}
    \label{eq:kkt2}
    \bar \Lambda^{\nu}_j G_j(\bar x) =0, \bar \Lambda^{\nu}_j \geq 0, G_j(\bar x) \geq 0, \quad j \in \mathcal{J},
\end{equation}
and, additionally:
\begin{equation}
    \label{eq:kkt3}
    r^1\bar \Lambda^1_j= \ldots = r^N \bar \Lambda^N_j, \quad j \in \mathcal{J}.
\end{equation}
We say that the normalized Nash equilibrium $\bar x$ is associated with $r$.
\end{definition}

It is clear that (\ref{eq:kkt1})-(\ref{eq:kkt2}) are Karush-Kuhn-Tucker optimality conditions for $\bar x^\nu$ to solve the parametric nonlinear programming problem $P_\nu(\bar x^{-\nu})$. We note that for stating Karush-Kuhn-Tucker condition some constraint qualification are usually needed.
The condition (\ref{eq:kkt3}) means that the multipliers for the same shared inequality constraint $G_j$ corresponding to different players $\nu \in \mathcal{N}$ are related, and this relation does not dependent on the index $j$. In particular, if $r^1=\ldots=r^N$, then the multipliers $\bar \Lambda^\nu_j$, $\nu \in \mathcal{N}$, coincide, even in case that the corresponding constraint $G_j$ is active. We point out that the notion of normalized Nash equiulibrium is invariant under a positive scaling of the parameters $r$. In what follows, we consider the parameters $r$ up to a positive scaling without mentioning this issue again and again.

The following GNEP-tailored constraint qualification is crucial for our analysis of normalized Nash equilibria.

\begin{definition}[GNEP-LICQ]
\label{def:gnep-licq}
We say that a GNEP feasible point $\bar x \in X$ satisfies the linear independence constraint qualification (GNEP-LICQ)  if the following system of linear equations 
\[
   \sum_{j \in \mathcal{J}^\nu_0(\bar x^\nu)}  \lambda^\nu_j D g_j^\nu(\bar x^\nu)+\sum_{j \in \mathcal{J}_0(\bar x)} \Lambda_j D_{x^\nu} G_j(\bar x) =0, \quad \nu \in \mathcal{N},
\]
admits only the trivial solution: 
\[
  \lambda^\nu_j = 0, j \in \mathcal{J}^\nu_0(\bar x^\nu), \nu \in \mathcal{N}, \quad  \Lambda_j=0, j \in \mathcal{J}_0(\bar x).
\]
\end{definition}

It turns out that GNEP-LICQ is a rather mild assumption.

\begin{theorem}[Genericity of GNEP-LICQ]
\label{thm:licq-gen}
Let $\mathcal{D}$  be the subset of GNEP defining functions $C^2\left(\mathbb{R}^n,\mathbb{R}^N\right)\times \prod_{\nu=1}^{N} C^2\left(\mathbb{R}^{n_\nu},\mathbb{R}^{\left|\mathcal{J}^\nu\right|}\right) \times C^2\left(\mathbb{R}^n,\mathbb{R}^{\left|\mathcal{J}\right|}\right)$ for which all feasible points satisfy GNEP-LICQ. Then, 
$\mathcal{D}$ is $C_s^2$-open and \mbox{-dense}.
\end{theorem}

\proof
For $(x,y) \in \mathbb{R}^{2n}$ we define the so-called structured $1$-jet extension:
\[
j^1(x,y)=\left(x^\nu, g^\nu_j(x^\nu), j \in \mathcal{J}^\nu, \nu  \in \mathcal{N}, y, G_j(y), j \in \mathcal{J}\right),
\]
which maps to $\mathbb{R}^{2n+\sum_{\nu=1}^{N} \left|\mathcal{J}^\nu\right|+\left|\mathcal{J}\right|}$ 
and corresponds to the functions
\[
  g=\left(g^\nu_j, j \in \mathcal{J}^\nu, \nu \in \mathcal{N}\right),\quad
   G=\left(G_j, j \in \mathcal{J} \right).
\]
Let us stratify the subset of the image space $ \mathbb{R}^{2n+\sum_{\nu=1}^{N} \left|\mathcal{J}^\nu\right|+\left|\mathcal{J}\right|}$ corresponding to the GNEP feasible set as follows:
\[
    A=\bigcup\limits_{\scriptsize \begin{array}{c}
         \mathcal{J}^\nu_1 \subset \mathcal{J}^\nu, \nu \in \mathcal{N}  \\
         \mathcal{J}_1 \subset \mathcal{J}
    \end{array}} A_{\mathcal{J}^1_1, \ldots, \mathcal{J}^N_1, \mathcal{J}_1},
\]
where we denote the variables by bold letters:
\[
  A_{\mathcal{J}^1_1, \ldots, \mathcal{J}^N_1, \mathcal{J}_1} = \left\{ (\mathbf{x}, \mathbf{g}, \mathbf{y}, \mathbf{G})  \,\left|\, 
  \begin{array}{l}
       \mathbf{x}^\nu = \mathbf{y}^\nu, \nu \in \mathcal{N},  \\
        \mathbf{g}^\nu_j=0, j \in \mathcal{J}^\nu \backslash \mathcal{J}^\nu_1, \\
        \mathbf{g}^\nu_j>0, j \in \mathcal{J}^\nu_1, \nu \in \mathcal{N}, \\
        \mathbf{G}_j=0, j \in \mathcal{J} \backslash \mathcal{J}_1, \\ \mathbf{G}^\nu_j>0, j \in \mathcal{J}_1
  \end{array}
    \right.\right\}.
\]
We show that GNEP-LICQ equivalently means that $j^1$ meets $A$ transversally.
Recall that $j^1$ meets the Whitney stratified set $A$ transversally if for all $(\bar x,\bar y)$ from the pre-image of $A$ under $j^1$ it holds:
\[
  D j^1(\bar x,\bar y)\left[\mathbb{R}^{2n}\right]+\mathcal{T}_{j^1(\bar x,\bar y)} A_{\mathcal{J}^1_1, \ldots, \mathcal{J}^N_1, \mathcal{J}_1}=\mathbb{R}^{2n+\sum_{\nu=1}^{N} \left|\mathcal{J}^\nu\right|+\left|\mathcal{J}\right|},
\]
where $A_{\mathcal{J}^1_1, \ldots, \mathcal{J}^N_1, \mathcal{J}_1}$ is the stratum of $A$ containing $j^1(\bar x,\bar y)$.
The differential $D j^1(\bar x,\bar y)\left[\mathbb{R}^{2n}\right]$ is spanned by the columns of the matrix
\[
B_1=\left(
\begin{array}{cccc}
e_i, i=1,\ldots, n_1 & &0 & 0\\ 
D g^1_j(\bar x^1), j \in\mathcal{J}^1 & & 0& 0 \\
 & \ddots & & \\ 
0 & &e_i, i=1,\ldots, n_N & 0\\ 
0 & & D g^N_j(\bar x^N), j \in \mathcal{J}^N& 0 \\
0 & &0 & e_i, i=1,\ldots, n\\ 
0 & & 0& D G_j(\bar y), j \in \mathcal{J} \\
\end{array}\right),
\]
The tangent space $\mathcal{T}_{j^1(\bar x,\bar y)} A_{\mathcal{J}^1_1, \ldots, \mathcal{J}^N_1, \mathcal{J}_1}$ is spanned by the columns of the matrix
\[
B_2=\left(
\begin{array}{cccc}
0 & &0 & 0\\ 
e_j, j \in\mathcal{J}^1_1 & & 0& 0 \\
 & \ddots & & \\ 
0 & &0 & 0\\ 
0 & & e_j, j \in \mathcal{J}^N_1& 0 \\
0 & &0 & 0\\ 
0 & & 0& e_j, j \in \mathcal{J}_1 \\
\end{array}\right)
\]
together with the columns of the matrix
\[
B_3=\left(
\begin{array}{ccc}
e_i, i=1,\ldots, n_1 & &0 \\ 
0 & & 0 \\
 & \ddots &  \\ 
0 & &e_i, i=1,\ldots, n_N \\ 
0 & & 0 \\
e_i, i=1,\ldots, n_1 & &0 \\ 
0 & & 0 \\
 & \ddots &  \\ 
0 & &e_i, i=1,\ldots, n_N \\ 
0 & & 0 \\
\end{array}\right).
\]
In order to show that $Dj^1(\bar x, \bar y)\left[\mathbb{R}^{2n}\right]$ and $\mathcal{T}_{j^1(\bar x,\bar y)} A_{\mathcal{J}^1_1, \ldots, \mathcal{J}^N_1, \mathcal{J}_1}$ sum up to the whole space $\mathbb{R}^{2n+\sum_{\nu=1}^{N} \left|\mathcal{J}^\nu\right|+\left|\mathcal{J}\right|}$, we determine the column rank of the composed matrix $B=\left(B_1, B_2, B_3\right)$.
Its column rank is equal to $2n+\sum_{\nu=1}^{N} \left|\mathcal{J}^\nu\right|+\left|\mathcal{J}\right|$ if and only if all its rows are linearly independent. This is exactly the case if just the following parts of the rows of $B$ are linearly independent:
\[
\left(
\begin{array}{ccc}
D g^1_j(\bar x^1), j \in \mathcal{J}^1 \backslash \mathcal{J}^1_1 & & 0  \\
 & \ddots &  \\ 
0 &  & D g^N_j(\bar x^N), j \in \mathcal{J}^N \backslash \mathcal{J}^N_1\\
D_{x^1} G_j(\bar y), j \in \mathcal{J} \backslash \mathcal{J}_1 & \ldots & D_{x^N} G_j(\bar y), j \in \mathcal{J} \backslash \mathcal{J}_1
\end{array}\right)
\]
Due to $j^1(\bar x,\bar y) \in A_{\mathcal{J}^1_1, \ldots, \mathcal{J}^N_1, \mathcal{J}_1}$, we have:
\[
\bar x = \bar y, \quad \mathcal{J}^\nu_0(\bar x)=\mathcal{J}^\nu \backslash \mathcal{J}^\nu_1, \nu \in \mathcal{N}, \quad \mathcal{J}_0(\bar x)=\mathcal{J} \backslash \mathcal{J}_1.
\]
By using this facts, transversality condition translates then into the linear independence of the rows of the matrix 
\[
\left(
\begin{array}{ccc}
D g^1_j(\bar x^1), j \in \mathcal{J}^1_0(\bar x) & & 0  \\
 & \ddots &  \\ 
0 &  & D g^N_j(\bar x^N), j \in  \mathcal{J}^N_0(\bar x)\\
D_{x^1} G_j(\bar x), j \in  \mathcal{J}_0(\bar x) & \ldots & D_{x^N} G_j(\bar x), j \in  \mathcal{J}_0(\bar x)
\end{array}\right).
\]
This is nothing else but GNEP-LICQ at $\bar x$.

Now, we apply the structured jet transversality theorem from \cite{guenzel:2008} to conclude the proof. Indeed, the latter says that for a given reduced jet extension and a given stratification the subset of functions, which meet the stratification transversally, is $C^2_s$-dense. For closed stratified sets it also gives that the mentioned subset of functions is $C^1_s$-open. Obviously, it is then also $C_s^2$-open. 
\qed

Let us briefly comment on how GNEP-LICQ is related to the standard LICQ from nonlinear programming.

\begin{remark}[Standard LICQ]
  GNEP-LICQ at $\bar x \in X$ equivalently requires the linear independence of the following vectors:
  \[
     D_x g_j^\nu(\bar x^\nu), j \in \mathcal{J}^\nu_0(\bar x^\nu), \nu \in \mathcal{N}, \quad D G_j(\bar x), j \in \mathcal{J}_0(\bar x),
  \]
where the derivatives of the individual constraints $g_j^\nu$'s are taken with respect to all variables in $x$ rather than just to $x^\nu$. Since the individual constraints do not depend on other player's variables, the genericity of GNEP-LICQ cannot be directly deduced from the standard results for nonlinear programming, cf. \cite{jongen:2000}. This comes from the fact that the derivatives of the individual constraints with respect to the other players' variables do trivially vanish, i.e.
\[
   D_{x^{-\nu}} g_j^\nu(\bar x^\nu) \equiv 0, j \in \mathcal{J}^\nu_0(\bar x^\nu).
\]
Hence, there is no freedom to perturb the whole gradients of $g_j^\nu$'s with respect to all $x$-variables, in order to eventually achieve their linear independence with the gradients of $G_j$'s. Nevertheless, Theorem \ref{thm:genericI} says that even in this restricted case a GNEP-tailored analogue of the standard LICQ generically holds. \qed
\end{remark}

It is straightforward to see that GNEP-LICQ implies GNEP-MFCQ as defined just below.

\begin{definition}[GNEP-MFCQ]
\label{def:gnep-mfcq}
We say that a GNEP feasible point $\bar x \in X$ satisfies the Mangasarian-Fromovitz constraint qualification (GNEP-MFCQ)  if there exists a vector $\xi=\left(\xi^1, \ldots, x^N\right) \in \mathbb{R}^n$, such that it holds: 
\[
   D g_j^\nu(\bar x^\nu) \xi^\nu > 0, j \in \mathcal{J}^\nu_0(\bar x^\nu), \nu \in \mathcal{N}, \quad
   D G_j(\bar x) \xi >0, j \in \mathcal{J}_0(\bar x).
\]
\end{definition}

Let us consider the subclass of GNEPs fulfilling the following standard convexity assumptions, \cite{rosen:1965}: 
\begin{itemize}
    \item[(C1)] $f^\nu(x^\nu,x^{-\nu})$ are convex in $x^\nu$ for each fixed $x^{-\nu}$, $\nu \in \mathcal{N}$,
    \item[(C2)] $g_j^\nu(x^\nu)$, $j \in \mathcal{J}^\nu$, $\nu \in \mathcal{N}$, and $G_j(x)$, $j \in \mathcal{J}$, are concave.    
\end{itemize}
Under (C2), the GNEP feasible set $X$ is convex. It is easy to see that in this case GNEP-MFCQ is equivalent to GNEP-Slater as defined below, cf. Theorem 3.2.75, \cite{stein:2018}.

\begin{definition}[GNEP-Slater]
\label{def:gnep-slater}
Let the convexity assumption (C2) be fulfilled.
We say that the GNEP feasible set $X$ satisfies the Slater constraint qualification (GNEP-Slater) if there exists $\widetilde x \in X$ with 
\[
  g_j^\nu(\widetilde x^\nu) > 0, j \in \mathcal{J}^\nu, \quad
   G_j(\widetilde x) >0, j \in \mathcal{J}.
\]
\end{definition}

Let us cite the main existence result on normalized Nash equilibria from \cite{rosen:1965}. We decided also to give its proof. Here, not only the decisive role of GNEP-tailored constraint qualifications becomes clear, but also the hint to studying the normalized Nash equilibria in the nonconvex setting is provided.    

\begin{theorem}[Existence of normalized Nash equilibrium, \cite{rosen:1965}]
\label{thm:ex-norm}
Let assumptions (C1) and (C2) be fulfilled, the GNEP feasible set $X$ be bounded and satisfy GNEP-Slater. Then, there exists a normalized Nash equilibrium associated with any parameters $r>0$.
\end{theorem}

\proof Let us define the following mapping for $(x,y) \in X \times X$:
\[
  \rho(x,y)= \sum_{\nu=1}^{N} r^\nu f^\nu(y^\nu, x^{-\nu}).
\]
Consider the point-to-set mapping for $x \in X$:
\[
  \Gamma(x) = \mbox{arg}\min_{y \in X} \rho(x,y).
\]
It follows from (C1) and (C2)
that $\Gamma$ is an upper semicontinuous mapping that maps each point of the convex,
compact set $X$ into a closed convex subset of $X$. Then, by the Kakutani fixed point
theorem \cite{kakutani:1941}, there exists a point $\bar x \in X$ such that $\bar x \in \Gamma(\bar x)$, i.e.
\[
  \bar x \in \mbox{arg}\min_{y \in X} \rho(\bar x,y).  
\]
Due to GNEP-Slater, Karush-Kuhn-Tucker optimality conditions can be stated for $\bar x$, i.e for all $\nu \in \mathcal{N}$ it holds: 
\begin{equation}
    \label{eq:h1}
r^\nu D_{x^\nu} f^\nu(\bar x) = \sum_{j \in \mathcal{J}^\nu} \bar \lambda^\nu_j D g_j^\nu(\bar x^\nu)+\sum_{j \in \mathcal{J}} \bar \Lambda_j D_{x^\nu} G_j(\bar x),
\end{equation}
\begin{equation}
    \label{eq:h2}
\bar \lambda^{\nu}_j g_j^\nu(\bar x^\nu) =0, \bar \lambda^{\nu}_j \geq 0, g_j^\nu(\bar x^\nu) \geq 0, \quad j \in \mathcal{J}^\nu,
\end{equation}
and, additionally,
\begin{equation}
    \label{eq:h3}
\bar \Lambda_j G_j(\bar x) =0, \bar \Lambda_j \geq 0, G_j(\bar x) \geq 0, \quad j \in \mathcal{J}.
\end{equation}
The formulae (\ref{eq:h1})-(\ref{eq:h3}) -- after dividing through $r^\nu$ -- become Karush-Kuhn-Tucker optimality conditions for $P_\nu(\bar x^{-\nu})$. Their sufficiency in convex programming provides that $\bar x^\nu$ solves $P_\nu(\bar x^{-\nu})$.
By setting
\[
       \bar \Lambda^\nu_j = \frac{\bar \Lambda_j}{r^\nu} \quad j \in \mathcal{J}, \nu \in \mathcal{N},
\]
the Nash equilibrium $\bar x$ is normalized.
\qed

Although the existence result from Theorem \ref{thm:ex-norm} is appealing, we emphasize that not all Nash equilibria are normalized, even if we allow in Definition \ref{def:norm-nash} to arbitrarily vary the parameters $r$. This obstacle for the notion of normalized Nash equilibrium has been already mentioned in Proposition 3.4, \cite{nabetani:2011}. The reasons for this are twofold:
\begin{itemize}
    \item[(a)]  Mangasarian-Fromovitz constraint qualification (MFCQ) may be violated at $\bar x^\nu$ for $P_\nu(\bar x^{-\nu})$. As a consequence, we cannot guarantee that Karush-Kuhn-Tucker optimality condition in (\ref{eq:kkt1}) would hold. Rather than that Fritz-John optimality condition happen to be valid.
    \item[(b)] Strict complementarity in (\ref{eq:kkt2}) need not to hold, i.e. some multipliers $\bar \Lambda^\nu_j$ corresponding to active shared constraints $G_j$ may vanish. Hence, if at least one of them does not vanish for the fixed $j$, the fulfillment of condition (\ref{eq:kkt3}) leads to a contradiction.
\end{itemize}
Let us illustrate the issues (a) and (b) by means of the following Examples (\ref{ex:issue1}) and (\ref{ex:issue2}), respectively. Note that the assumptions of Theorem \ref{thm:ex-norm} hold here, i.e. (C1) and (C2) are fulfilled, the GNEP feasible set $X$ is bounded and satisfies GNEP-Slater. However, these examples exhibit Nash equilibria, which are not normalized, and this phenomenon is stable with respect to arbitrarily small perturbations of the defining functions. 

\begin{example}[Failure of KKT condition, cf. \cite{dorsch:2013}]
  \label{ex:issue1}
  Let $\mathcal{N} = \{1, 2\}$ and GNEP
be given by
\[
   f^1\left(x^1_1,x^1_2,x^2\right)=-x^1_1, \quad f^2\left(x^1_1,x^1_2,x^2\right)=-x^2,
\]
\[
   G_1\left(x^1_1,x^1_2,x^2\right)=1-\left(x^1_1-x^2\right)^2-\left(x^1_2-\left(1-2x^2\right)\right)^2,
\]
\[
    G_2\left(x^1_1,x^1_2,x^2\right)=1-\left(x^1_1\right)^2-\left(x^2+1\right)^2.
\]
Let us consider its Nash equilibrium $\bar x =\left(0,0,0\right)$. The Nash equilibrium $\bar x$ is not normalized, since $\left(\bar x^1_1, \bar x^1_2\right)=(0,0)$ is not a Karush-Kuhn-Tucker point for $P_1\left(\bar x^{-1}\right)=P_1(0)$, but just a Fritz-John point. Hence, (\ref{eq:kkt1}) cannot hold for any parameters $r$. This matter is due to the violation of MFCQ for the first player's optimization problem at the origin. 
\qed
\end{example}

\begin{example}[Violation of strict complementarity]
  \label{ex:issue2}
Let $\mathcal{N} = \{1, 2\}$ and GNEP
be given by
\[
   f^1\left(x^1,x^2\right)= \left(x^1-x^2\right)^2, \quad f^2\left(x^1,x^2\right)=-x^2,
\]
\[
   G_1(x^1,x^2)=x^1-2x^2, \quad G_2(x^1,x^2)=1-x^1, \quad G_3(x^1,x^2)=x^2+1.
\]
Let us consider its Nash equilibrium $\bar x =\left(0,0\right)$. To see that $\bar x$ is not normalized, we compute the multipliers of $\bar x^1=0$ for $P_1\left(\bar x^{-1}\right)=P_1(0)$ and of $\bar x^2=0$ for $P_2\left(\bar x^{-2}\right)=P_2(0)$ corresponding to the active constraints from $\mathcal{J}_0(\bar x)=\left\{1\right\}$, respectively:
\[
   \bar \Lambda^1_1 =0, \quad 
   \bar \Lambda^2_1 = \frac{1}{2}.
\]
Hence, (\ref{eq:kkt2}) cannot hold for any parameters $r$. This matter is due to the violation of strict complementarity for the first player's optimization problem at the origin.
\qed  
\end{example}


At the end of this section we put the notion of normalized Nash equilibrium into the context of our previous study \cite{dorsch:2013} on the structure of Nash equilibria for GNEP.

\begin{remark}[Additional equations]
  \label{rem:study}
  In \cite{dorsch:2013}, the structure of Nash equilibria for GNEP in presence of shared constraints has been addressed. It turns out that the set of Nash equilibria -- considered as Fritz-John points  together with the corresponding multipliers -- constitutes generically a Lipschitz manifold. Its dimension locally at a Nash equilibrium $\bar x \in X$ is $(N-1)\left|\mathcal{J}_0(\bar x)\right|$, see Corollary 2.4, \cite{dorsch:2013}. This means that the concatenated system of Fritz-John conditions for $\bar x^\nu$ as solutions of $P_\nu(\bar x^{-\nu})$, $\nu \in \mathcal{N}$, has more variables than equations. This "deficiency" explains the geometrical complexity of the set of Nash equilibria, such as its non-uniqueness, the appearance of kinks and boundary points, \cite{dorsch:2013a}. From this point of view, the notion of normalized Nash equilibrium just introduces additional equations (\ref{eq:kkt3}) into the concatenated system of Karush-Kuhn-Tucker conditions. Note that their number in (\ref{eq:kkt3}), after omitting trivial ones corresponding to inactive shared inequality constraints, amounts to exactly $(N-1)\left|\mathcal{J}_0(\bar x)\right|$. Since now the numbers of variables and equations in (\ref{eq:kkt1})-(\ref{eq:kkt3}) coincide, we expect normalized Nash equilibria to become locally unique, at least for a generic GNEP. \qed
\end{remark}

\section{Local uniqueness}
\label{sec:lu}

From the proof of Theorem \ref{thm:ex-unique} we see how to alternatively describe normalized Nash equilibria.
Equivalently to fulfilling (\ref{eq:kkt1})-(\ref{eq:kkt3}), a normalized Nash equilibrium $\bar x \in X$ can be defined by requiring the following conditions to hold for all $\nu \in \mathcal{N}$:
\begin{equation}
    \label{eq:kkt1h}
    r^\nu D_{x^\nu} f^\nu(\bar x) = \sum_{j \in \mathcal{J}^\nu} \bar \lambda^\nu_j D g_j^\nu(\bar x^\nu)+\sum_{j \in \mathcal{J}} \bar \Lambda_j D_{x^\nu} G_j(\bar x),
\end{equation}
\begin{equation}
    \label{eq:kkt2ah}
    \bar \lambda^{\nu}_j g_j^\nu(\bar x^\nu) =0, \bar \lambda^{\nu}_j \geq 0, g_j^\nu(\bar x^\nu) \geq 0, \quad j \in \mathcal{J}^\nu,
\end{equation}
and, additionally,
\begin{equation}
    \label{eq:kkt2h}
    \bar \Lambda_j G_j(\bar x) =0, \bar \Lambda_j \geq 0, G_j(\bar x) \geq 0, \quad j \in \mathcal{J}.
\end{equation}
This is possible since, by virtue of (\ref{eq:kkt3}), we may set as the new multipliers corresponding to the shared constraints:
\[
    \bar \Lambda_j = r^1 \bar \Lambda^1_j = \ldots = r^N \bar \Lambda^N_j, \quad j \in \mathcal{J}.
\]
The multipliers $\bar \lambda^{\nu}_j$,$j \in \mathcal{J}^\nu(\bar x^\nu)$, which correspond to the individual constraints, are just scaled  by means of $r^\nu$ accordingly. By formalizing this discussion, we arrive at the notion of normalized Karush-Kuhn-Tucker point for GNEP.

\begin{definition}[Normalized KKT-point] \label{def:norm-kkt}
A GNEP feasible point $\bar x \in X$ is called normalized Karush-Kuhn-Tucker (KKT) point in association with the parameters $r>0$, if it fulfills conditions (\ref{eq:kkt1h})-(\ref{eq:kkt2h}) with the corresponding multipliers.
\end{definition}


If GNEP-LICQ holds at a normalized KKT-point $\bar x \in X$, then the multipliers
\[ 
  \bar \lambda = \left(\bar \lambda^{\nu}_j, j \in \mathcal{J}^\nu_0(\bar x^\nu), \nu \in \mathcal{N}\right), \quad \bar \Lambda = \left(\bar \Lambda_j, j \in \mathcal{J}_0(\bar x) \right).
\]
are uniquely determined. Here, we skip zero multipliers corresponding to the inactive constraints. It is convenient to also define the Lagrange function for every $\nu \in \mathcal{N}$:
    \[
     L^\nu(x,\lambda, \Lambda,r)=r^\nu f^\nu(x) - \sum_{j \in \mathcal{J}^\nu_0(\bar x^\nu)} \lambda^\nu_j g_j^\nu(x^\nu)-\sum_{j \in \mathcal{J}_0(\bar x)} \Lambda_j G_j(x).
    \]
Conditions in (\ref{eq:kkt1h}) read then as $D_{x^\nu}L^\nu(\bar x,\bar \lambda, \bar \Lambda,r) =0$, $\nu \in \mathcal{N}$. Equivalently, we write $\mathcal{G}_L(\bar x,\bar \lambda, \bar \Lambda,r)=0$ for the latter, where the so-called pseudogradient mapping with respect to $x$-variables is defined as follows:
\[
   \mathcal{G}_L(x, \lambda, \Lambda,r)= \left(D_{x^\nu}L^\nu(x,\lambda, \Lambda,r), \nu \in \mathcal{N}\right)^T.
\]
Its Jacobian $D_x \mathcal{G}_L(x,\lambda, \Lambda,r)$ is an $(n \times n)$-matrix.


Further, we set
\[
  X_0(\bar x) =\left\{
x \in \mathbb{R}^n\,\left\vert\,\begin{array}{l}
     g_j^\nu(x^\nu) =0, j \in\mathcal{J}^\nu_0(\bar x^\nu), \nu \in \mathcal{N} \\
G_j(x)=0, j \in\mathcal{J}_0(\bar x)
\end{array} 
\right.\right\}.
\]
Obviously, locally at $\bar x$ we have $X_0(\bar x) \subset X$. In case that GNEP-LICQ holds at $\bar x$, the set $X_0(\bar x)$ is locally an $C^2$-manifold of dimension $n- \sum_{\nu=1}^{N} \left|J^\nu_0(\bar x^\nu)\right|-\left|J_0(\bar x)\right|$. 
The tangent space of $X_0(\bar x)$ at $\bar x$ is thus given by
\[
    \mathcal{T}_{\bar x} X_0(\bar x)=\left\{
\xi \in \mathbb{R}^n\,\left\vert\, \begin{array}{l} Dg_j^\nu(x^\nu) \xi^\nu = 0, j \in\mathcal{J}^\nu_0(\bar x^\nu), \nu \in \mathcal{N} \\
DG_j(x) \xi=0, j \in\mathcal{J}_0(\bar x)
\end{array}
\right.\right\}.
\]

We define the concept of nondegeneracy for normalized KKT-points. As we shall see later, this is a sufficient condition to guarantee that normalized KKT-points are locally unique.

\begin{definition}[Nondegenerate normalized KKT-point]
\label{def:nd-m-st}
A normalized KKT-point $\bar x \in X$ associated with $r>0$ is called nondegenerate if
\begin{itemize}
    \item[] ND1: GNEP-LICQ is fulfilled at $\bar x$,
    \item[] ND2: $\bar \lambda^\nu_j>0$ for all $j\in \mathcal{J}^\nu_0(\bar x^\nu)$ and $\bar \Lambda_j>0$ for all $j\in \mathcal{J}_0(\bar x)$,
    \item[] ND3: the matrix $D_x \mathcal{G}_L(\bar x, \bar \lambda, \bar \Lambda,r)\restriction_{\mathcal{T}_{\bar x} X_0(\bar x)}$ is nonsingular, i.e. the matrix $V^T D_x \mathcal{G}_L(\bar x, \bar \lambda, \bar \Lambda,r) V$ is nonsingular, where $V$ is some matrix whose columns form a basis of the tangent space $\mathcal{T}_{\bar x} X_0(\bar x)$.      
    \end{itemize}
\end{definition}

Nondegeneracy of a normalized Nash equilibrium $\bar x \in X$ does not in general prevent from degeneracies at the individual level. The latter refers to the possible violation of the standard linear independence constraint qualification (LICQ), strict complementarity (SC), or second-order sufficient condition (SOSC) at $\bar x^\nu$ as solution of $P_\nu(\bar x^{-\nu})$ for some $\nu \in \mathcal{N}$. This issue is illustrated in Example \ref{ex:deg}.

\begin{example}[Individual degeneracies]
\label{ex:deg}
Let $\mathcal{N} = \{1, 2\}$ and GNEP
be given by
\[
   f^1(x^1,x^2)=-x^1, \quad f^2(x^1,x^2)=-x^2,
\]
\[
   G_1(x^1,x^2)=1-x^1-x^2, \quad
    G_2(x^1,x^2)=x^1-x^2, \quad
     G_3(x^1,x^2)=x^2.
\]
Let us consider its Nash equilibrium $\bar x =\left(\frac{1}{2},\frac{1}{2}\right)$. It is normalized in association with any parameters $r=(r^1,r^2)$, fulfilling $\frac{r^1}{r^2}<1$. The corresponding multipliers for the active constraints from $\mathcal{J}_0(\bar x)=\left\{1,2\right\}$ are 
\[
  \bar \Lambda^1=\frac{1}{2}r^1+\frac{1}{2}r^2, \quad
  \bar \Lambda^2=-\frac{1}{2}r^1+\frac{1}{2}r^2.
\]
Obviously, GNEP-LICQ is fulfilled at $\bar x$, i.e. also ND1. Due to the choice of the parameters $r$, ND2 holds at $\bar x$ as well. Since $\mathcal{T}_{\bar x} X_0(\bar x)=\{0\}$, ND3 is trivially satisfied. Overall, the normalized Nash equilibrium $\bar x$ is nondegenerate due to Definition \ref{def:nd-m-st}. However, the standard LICQ is violated at the solution $\bar x^1=\frac{1}{2}$ for the first player's optimization problem $P_1(\bar x^{-1})=P_1\left(\frac{1}{2}\right)$ and 
at the solution $\bar x^2=\frac{1}{2}$ for the second player's optimization problem $P_2(\bar x^{-2})=P_2\left(\frac{1}{2}\right)$. Even worse, the strict complementarity need not to hold for both individual optimization problems $P_1(\bar x^{-1})$ and $P_2(\bar x^{-2})$ at $\bar x^1$ and $\bar x^2$, respectively. Here, each of the usual KKT conditions may be separately satisfied with some vanishing multipliers corresponding to the active shared constraints. We emphasize that the phenomenon presented here is stable with respect to arbitrarily small perturbations of the defining functions.
\qed
\end{example}

It turns out that nondegeneracy is a generic property of normalized KKT-points for GNEP associated with particular parameters $r$.

\begin{theorem}[Genericity of nondegeneracy]
\label{thm:genericI}
Let $\mathcal{D}$  be the subset of GNEP defining functions $C^2\left(\mathbb{R}^n,\mathbb{R}^N\right)\times \prod_{\nu=1}^{N} C^2\left(\mathbb{R}^{n_\nu},\mathbb{R}^{\left|\mathcal{J}^\nu\right|}\right) \times C^2\left(\mathbb{R}^n,\mathbb{R}^{\left|\mathcal{J}\right|}\right)$ for which each 
normalized KKT-point associated with any, but fixed, parameters $r>0$ is nondegenerate. Then, $\mathcal{F}$ is $C^2_s$-open and -dense.
\end{theorem}

\proof  Let us fix index subsets $\mathcal{J}_0^\nu \subset \mathcal{J}^\nu$, $\nu \in \mathcal{N}$, $\mathcal{J}_0 \subset \mathcal{J}$ of active individual and shared inequality constraints, index subsets $\mathcal{K}_0^\nu \subset \mathcal{J}^\nu_0$, $\nu \in \mathcal{N}$, $\mathcal{K}_0 \subset \mathcal{J}_0$ of these active inequality constraints, and a number $\ell \in\{0,1, \ldots, n\}$ standing for the rank.
For this choice we consider the set  $M$ of $x \in \mathbb{R}^n$ such that the following conditions are satisfied:
\begin{itemize}
    \item[] (m1) $g^\nu_j(x^\nu)=0$, $j \in \mathcal{J}^\nu_0$, $g^\nu_j(x^\nu)>0$, $j \in \mathcal{J}^\nu \backslash \mathcal{J}^\nu_0$, $\nu \in \mathcal{N}$, and 
    $G_j(x)=0$, $j \in \mathcal{J}_0$,  
    $G_j(x)>0$, $j \in \mathcal{J} \backslash \mathcal{J}_0$,
    \item[] (m2) the vector 
    \[
       \left( \begin{array}{c}
            r^1 D_{x^1}^T f^1(x)  \\
            \vdots \\
            r^\nu D_{x^\nu}^T f^\nu(x)\\
            \vdots\\
            r^N D_{x^N}^T f^N(x)  
       \end{array}\right) 
    \]
    is spanned by the vectors
       \[
       \left( \begin{array}{c}
            0  \\
            \vdots \\
            D_{x^\nu}^T g^\nu_j(x^\nu)\\
            \vdots \\
            0 
       \end{array}\right), j \in \mathcal{J}_0^\nu \backslash \mathcal{K}_0^\nu, \nu \in \mathcal{N},\quad
         \left( \begin{array}{c}
            D_{x^1}^T G_j(x)  \\
            \vdots \\
            D_{x^\nu}^T G_j(x)\\
            \vdots\\
            D_{x^N}^T G_j(x) 
       \end{array}\right), j \in \mathcal{J}_0 \backslash \mathcal{K}_0,
    \]
    \item[] (m3) the matrix $D_x \mathcal{G}_L(x,\lambda, \Lambda,r)\restriction_{\mathcal{T}_{x} X_0(x)}$ has rank $\ell$.
\end{itemize}
Note that (m1) models feasibility and explicitly refers to active inequality constraints. Conditions (m2) and (m3) describe possible violation of ND2 and ND3, respectively.
Now, it suffices to show that $M$ is generically empty whenever one of the sets $\mathcal{K}_0^\nu$, $\nu \in \mathcal{N}$, or $\mathcal{K}_0$ is nonempty or the rank $\ell$ in (m3) is not full, i.\,e. $\ell < \mbox{dim}\left(\mathcal{T}_{x} X_0(x)\right)$.
In fact, the available degrees of freedom
of the variables involved in each $M$ are $n$. The loss of freedom caused by (m1) is $\sum_{\nu=1}^{N} \left|\mathcal{J}_0^\nu\right|+\left|\mathcal{J}_0\right|$.
Due to Theorem \ref{thm:licq-gen}, GNEP-LICQ holds generically at any GNEP feasible point $x$, i.\,e. (ND1) is fulfilled. In particular, the vectors forming the span in (m2) are linearly independent. Suppose that the sets $\mathcal{K}_0^\nu$, $\nu \in \mathcal{N}$, and $\mathcal{K}_0$ are empty, then (m2) causes a loss of freedom of $n-\sum_{\nu=1}^{N} \left|\mathcal{J}_0^\nu\right|-\left|\mathcal{J}_0\right|$. Hence, the total loss of freedom is $n$.
We conclude that a
further degeneracy, i.\,e. $\mathcal{K}_0^\nu \not = \emptyset$, $\nu \in \mathcal{N}$, $\mathcal{K}_0 \not = \emptyset$ or $\ell < \mbox{dim}\left(\mathcal{T}_{x} X_0(x)\right)$, would imply that the total available degrees of freedom $n$ are exceeded. By virtue of the jet
transversality theorem from \cite{jongen:2000}, generically the sets $M$ must be empty.
For the openness result, we argue in a standard way. Locally, normalized KKT-points can be written
via stable equations. Then, the implicit function theorem for Banach spaces can be applied to
follow normalized KKT-points with respect to (local) $C^2$-perturbations of defining functions. Finally,
a standard globalization procedure exploiting the specific properties of the strong $C^2_s$-topology can be used to construct a (global) $C^2_s$-neighborhood of problem data for which the nondegeneracy
property is stable, cf. \cite{jongen:2000}. \qed



Let us illustrate the assertion of Theorem \ref{thm:genericI} 
by means of Example \ref{ex:per}.

\begin{example}[Perturbation]
\label{ex:per}
Let $\mathcal{N} = \{1, 2\}$ and GNEP
be given again by
\[
   f^1(x^1,x^2)=-x^1, \quad f^2(x^1,x^2)=-x^2,
\]
\[
   G_1(x^1,x^2)=1-x^1-x^2, \quad
    G_2(x^1,x^2)=x^1-x^2, \quad
     G_3(x^1,x^2)=x^2.
\]
Let us consider its Nash equilibria $\bar x(t)=(1-t,t)$ with $t \in \left(0,\frac{1}{2}\right)$. They are normalized in association with the parameters $r=(r^1,r^2)$, fulfilling $\frac{r^1}{r^2}=1$. 
The corresponding multiplier for the active constraint from $\mathcal{J}_0(\bar x(t))=\left\{1\right\}$ is 
\[
  \bar \Lambda^1=r^1.
\]
Obviously, GNEP-LICQ is fulfilled at $\bar x(t)$, i.e. also ND1. The unique multiplier $\bar \Lambda^1$ is positive, hence, ND2 holds at $\bar x(t)$ as well. However, $\bar x(t)$ is degenerate, since ND3 is violated. Theorem \ref{thm:genericI} suggests that it is possible to perform an arbitrarily small $C^2$-perturbation of the defining functions of this GNEP, in order to guarantee that all normalized Nash equilibria associated with the parameters $r$, fulfilling $\frac{r^1}{r^2}=1$, are nondegenerate. Here, it is sufficient to perturb just the objective functions of GNEP as follows:
\[
   \widetilde f^1(x^1,x^2)=-x^1+\frac{\varepsilon}{2} x^1 \cdot x^1, \quad \widetilde f^2(x^1,x^2)=-x^2+\varepsilon x^2 \cdot x^2,
\]
where $\varepsilon > 0$ is taken sufficiently small.
It is technical, but not hard to see that $\bar x= \left(\frac{2}{3}, \frac{1}{3}\right)$ is the only normalized Nash equilibrium of the perturbed GNEP associated to the parameters $r=(r^1,r^2)$, fulfilling $\frac{r^1}{r^2}=1$. The corresponding multiplier for the active constraint from $\mathcal{J}_0(\bar x)=\left\{1\right\}$ is 

\[
  \bar \Lambda^1=\left(1-\frac{2}{3}\varepsilon\right)r^1.
\]
As previously, GNEP-LICQ is fulfilled at $\bar x$, i.e. also ND1. For a sufficiently small $\varepsilon >0$, ND2 holds at $\bar x$ as well. Finally, we have:
\[
  \mathcal{T}_{\bar x} X_0(\bar x) = \left\{(\xi^1,\xi^2)^T \,\left|\, -\xi^1-\xi^2=0\right.\right\}
\]
and
\[
  D_x \mathcal{G}_L(\bar x, \bar \Lambda, r)= \left(
  \begin{array}{cc}
      \varepsilon r^1 & 0 \\
      0 & 2\varepsilon r^2
  \end{array}
  \right).
\] 
Hence, ND3 is also satisfied, since
for the perturbed GNEP it holds:
\[
 \begin{array}{rcl}
      \displaystyle
  D_x \mathcal{G}_L(\bar x, \bar \Lambda, r)\restriction_{\mathcal{T}_{\bar x} X_0(\bar x)} &=& \displaystyle
  \left(
  \begin{array}{cc}
      1 & -1 
  \end{array}
  \right) \left(
  \begin{array}{cc}
      \varepsilon r^1 & 0 \\
      0 & 2\varepsilon r^2
  \end{array}
  \right) \left(
  \begin{array}{c}
      1 \\ -1 
  \end{array}
  \right) \\ \\ &=& \displaystyle \varepsilon r^1 + 2\varepsilon r^2 \not =0.
 \end{array}
\] 
Overall, the normalized Nash equilibrium $\bar x= \left(\frac{2}{3}, \frac{1}{3}\right)$ is nondegenerate. \qed 
\end{example}

It turns out that nondegenerate normalized KKT-points (in particular, normalized Nash equilibria) are locally unique.

\begin{theorem}[Nondegeneracy and local uniqueness]
\label{thm:unique}
Nondegenerate normalized KKT-points (in particular, normalized Nash equilibria) associated with any, but fixed, parameters $r>0$ are locally unique. 
\end{theorem}

\proof Let $\bar x \in X$ be a nondegenerate normalized KKT-point associated with parameters $r^\nu >0$, $\nu \in \mathcal{N}$. We denote its unique multipliers as
\[ 
  \bar \lambda = \left(\bar \lambda^{\nu}_j, j \in \mathcal{J}^\nu_0(\bar x^\nu), \nu \in \mathcal{N}\right), \quad \bar \Lambda = \left(\bar \Lambda_j, j \in \mathcal{J}_0(\bar x) \right).
\]
 Let us assume that there exists a sequence of normalized KKT-points $x(k)$, $k=1,2,\ldots$, associated with the fixed set of parameters $r^\nu >0$, $\nu \in \mathcal{N}$, such that $x(k) \rightarrow \bar x$ if $k \rightarrow \infty$. Due to continuity reasons, for sufficiently large $k$ it holds -- after taking a subsequence if needed: 
\[
   \mathcal{J}^\nu_0(x^\nu(k)) \subset \mathcal{J}^\nu_0(\bar x^\nu), \nu \in \mathcal{N}, 
   \quad \mathcal{J}_0(x(k)) \subset \mathcal{J}_0(\bar x).
\]
Thus, we may denote the unique multipliers of $x(k)$ as 
\[ 
  \lambda(k) = \left(\lambda^{\nu}_j(k), j \in \mathcal{J}^\nu_0(\bar x^\nu), \nu \in \mathcal{N}\right), \quad \Lambda(k) = \left( \Lambda_j(k), j \in \mathcal{J}_0(\bar x) \right).
\]
Let us show that the sequence of multipliers $(\lambda(k), \Lambda(k))$, $k \in \mathbb{N}$, is bounded. For that, we write down (\ref{eq:kkt1h}) as follows:

 \[
       \begin{array}{rcl}
       \left( \begin{array}{c}
            r^1 D_{x^1}^T f^1(x(k))  \\
            \vdots \\
            r^\nu D_{x^\nu}^T f^\nu(x(k))\\
            \vdots\\
            r^N D_{x^N}^T f^N(x(k))  
       \end{array}\right) &=&  
        \displaystyle
       \sum_{\nu \in \mathcal{N}}\sum_{j \in \mathcal{J}^\nu_0(\bar x^\nu)} \lambda^{\nu}_j(k)\left( \begin{array}{c}
            0  \\
            \vdots \\
            D_{x^\nu}^T g^\nu_j(x^\nu(k))\\
            \vdots \\
            0 
       \end{array}\right) \\ \\
       &&+ \displaystyle \sum_{j \in \mathcal{J}_0(\bar x)} \Lambda_j(k) \left( \begin{array}{c}
            D_{x^1}^T G_j(x(k))  \\
            \vdots \\
            D_{x^\nu}^T G_j(x(k))\\
            \vdots\\
            D_{x^N}^T G_j(x(k)) 
       \end{array}\right).
       \end{array}
    \]
We rewrite these system of equations by using the pseudogradient notation:
\begin{equation}
    \label{eq:kkt-short}
 \mathcal{G}_{f}(x(k),r) = \left(D^T g(x(k)), D^T G(x(k))\right) \cdot \left( \begin{array}{c}
        \lambda(k) \\ \Lambda(k)
   \end{array}\right).
   \end{equation}
Here, we put with some abuse of notation:
\[
  \quad \mathcal{G}_f(x,r) =\left(r^\nu D_{x^\nu}^T f^\nu(x), \nu \in \mathcal{N}\right),
\]
and 
\[
g(x)=\left(g^\nu_j(x^\nu), j \in \mathcal{J}^\nu_0(\bar x^\nu), \nu \in \mathcal{N}\right),\quad
   G(x)=\left(G_j(x), j \in \mathcal{J}_0(\bar x) \right).
\]
By applying the Euclidean norm to (\ref{eq:kkt-short}), we obtain:
\[
   \left\|\mathcal{G}_{f}(x(k),r) \right\| \geq \min_{\left\|(\lambda, \Lambda)\right\|_2=1}\left\| \left(D^T g(x(k)), D^T G(x(k))\right) \left( \begin{array}{c}
        \lambda \\ \Lambda
   \end{array}\right) \right\| \cdot \left\| \left( \begin{array}{c}
        \lambda(k) \\ \Lambda(k)
   \end{array}\right) \right\|.
\]
Due to GNEP-LICQ at $\bar x$ from (ND1), the matrix $\left(D^T g(\bar x), D^T G(\bar x)\right)$ has full column rank. Hence, since $x(k) \rightarrow \bar x$ if $k \rightarrow \infty$, the expressions
\[
   \min_{\left\|(\lambda, \Lambda)\right\|_2=1}\left\| \left(D^T g(x(k)), D^T G(x(k))\right) \left( \begin{array}{c}
        \lambda \\ \Lambda
   \end{array}\right) \right\|
\]
are bounded away from zero uniformly for sufficiently large $k$. Obviously, the convergent sequence $\left\|\mathcal{G}_{f}(x(k),r) \right\|$ is also bounded in $k$. The boundedness of $(\lambda(k), \Lambda(k))$, $k \in \mathbb{N}$, thus easily follows.

Now, we may assume without loss of generality -- by considering a subsequence if needed -- that the sequence $(\lambda(k), \Lambda(k))$, $k \in \mathbb{N}$, converges. Due to GNEP-LICQ at $\bar x$, it must then hold:
\[
(\lambda(k), \Lambda(k)) \rightarrow (\bar \lambda, \bar \Lambda) \mbox{ if } k \rightarrow \infty.
\]
In particular, we have for sufficiently large $k$, due to (ND2):
\[
\lambda^\nu_j(k) > 0, j \in \mathcal{J}^\nu_0(\bar x^\nu), \quad \Lambda_j(k) > 0, j \in \mathcal{J}_0(\bar x).
\]
Consequently, we obtain:
\[
   \mathcal{J}^\nu_0(x^\nu(k)) = \mathcal{J}^\nu_0(\bar x^\nu), \nu \in \mathcal{N}, 
   \quad \mathcal{J}_0(x(k)) = \mathcal{J}_0(\bar x).
\]
In order to come to a contradiction, let us consider the following mapping from $\mathbb{R}^{n+ \sum_{\nu=1}^{N} \left|\mathcal{J}^\nu_0(\bar x^\nu)\right|+\left|\mathcal{J}_0(\bar x)\right|}$ into itself:
\[
   \mathcal{F}(x,\lambda,\Lambda) =
   \left(
   \begin{array}{c}
        \displaystyle D_{x^\nu}^T L^\nu(x,\lambda, \Lambda,r), \nu \in \mathcal{N} \\
         g_j^\nu(x^\nu), j \in \mathcal{J}^\nu_0(\bar x^\nu),\nu \in \mathcal{N} \\ 
         G_j(x), j \in \mathcal{J}_0(\bar x)
   \end{array}
   \right),
\]
where the Lagrange functions for every $\nu \in \mathcal{N}$ are
    \[
     L^\nu(x,\lambda, \Lambda,r)=r^\nu f^\nu(x) - \sum_{j \in \mathcal{J}^\nu_0(\bar x^\nu)} \lambda^\nu_j g_j^\nu(x^\nu)-\sum_{j \in \mathcal{J}_0(\bar x)} \Lambda_j G_j(x).
    \]
Equivalently, we write for short:
\[
    \mathcal{F}(x,\lambda,\Lambda) =
   \left(
   \begin{array}{c}
        \displaystyle \mathcal{G}_L(x,\lambda,\Lambda,r)  \\
         g(x) \\ 
         G(x)
   \end{array}
   \right).
\]
Its Jacobian evaluated at $(\bar x,\bar \lambda, \bar \Lambda)$ is
\[
   D \mathcal{F}(\bar x,\bar \lambda, \bar \Lambda) = 
   \left(
      \begin{array}{ccc}
         D_x \mathcal{G}_L(\bar x,\bar \lambda,\bar \Lambda,r) & -D^Tg(\bar x) & -D^TG(\bar x) \\
         Dg(\bar x)  & 0 & 0 \\
         DG(\bar x)  & 0 & 0 \\
      \end{array}
   \right).
\]
In presence of GNEP-LICQ, the matrix $D \mathcal{F}(\bar x,\bar \lambda, \bar \Lambda)$ is nonsingular if and only if $D_x \mathcal{G}_L(\bar x,\bar \lambda,\bar \Lambda,r)\restriction_{\mathcal{T}_{\bar x} X_0(\bar x)}$ is nonsingular, cf. Theorem 2.3.2, \cite{jongen:2004}. The nonsingularity of $D_x \mathcal{G}_L(\bar x,\bar \lambda,\bar \Lambda,r)\restriction_{\mathcal{T}_{\bar x} X_0(\bar x)}$ holds due to (ND3). The inherited nonsingularity of the Jacobian $D \mathcal{F}(\bar x,\bar \lambda, \bar \Lambda)$ allows to apply the inverse function theorem at a zero $(\bar x,\bar \lambda, \bar \Lambda)$ of $\mathcal{F}$. In particular, $(\bar x,\bar \lambda, \bar \Lambda)$ is the locally unique solution of the system of equations $\mathcal{F}(x,\lambda, \Lambda)=0$. However, for sufficiently large $k$ we also have $\mathcal{F}(x(k),\lambda(k), \Lambda(k))=0$, a contradiction. \qed

The direct application of Theorems \ref{thm:genericI} and \ref{thm:unique} provides the following result.

\begin{corollary}[Genericity and local uniqueness]
\label{cor:unique}
Generically, all normalized KKT-points (in particular, all normalized Nash equilibria) associated with any, but fixed, parameters $r>0$, are locally unique. 
\end{corollary}

Now, we compare our results on the local uniqueness of normalized Nash equilibria with those due to the seminal paper \cite{rosen:1965}. 
There, the following assumption -- together with (C1) and (C2) -- is made for fixed parameters $r>0$:
\begin{itemize}
    \item[(C3)] the Jacobian $D_x \mathcal{G}_{f}(x,r)$ is positive definite for all $x \in X$.
\end{itemize}
Here, the so-called pseudogradient mapping with respect to $x$-variables is defined as follows:
\[
     \mathcal{G}_{f}(x,r)=\left(r^\nu D_{x^\nu} f^\nu(x), \nu \in \mathcal{N}\right)^T.
\]
Its Jacobian $D_x \mathcal{G}_f(x,r)$ is an $(n \times n)$-matrix.
We recall that the possibly nonsymmetric matrix $D_x \mathcal{G}_{f}(x,r)$ is positive definite if and only if its symmetric part $\frac{1}{2}\left(D_x \mathcal{G}_{f}(x,r)+D^T_x \mathcal{G}_{f}(x,r)\right)$ is positive definite. Condition (C3) is crucial for guaranteeing uniqueness of normalized Nash equilibria in the convex setting.

\begin{theorem}[Uniqueness of normalized Nash equilibrium, \cite{rosen:1965}]
\label{thm:ex-unique}
Let assumptions (C1), (C2) be fulfilled, the GNEP feasible set $X$ be bounded and satisfy GNEP-Slater. If (C3) holds for any, but fixed, parameters $r>0$, then there is a unique normalized Nash equilibrium associated with $r$.
\end{theorem}

It turns out that the proposed notion of nondegeneracy from Definition \ref{def:nd-m-st} and condition (C3) are of independent interest, if studying (local) uniqueness of normalized Nash equilibria. As next Example \ref{ex:compar} shows, they are not reducible to each other, even if the convexity assumptions (C1) and (C2) are additionally met.  

\begin{example}[Nondegeneracy vs. (C3)]
\label{ex:compar}
Let $\mathcal{N} = \{1, 2\}$ and GNEP
be given by
\[
   f^1(x^1,x^2)=-x^1+x^1\cdot x^2, \quad f^2(x^1,x^2)=-x^2+\frac{1}{2}x^1\cdot x^2,
\]
\[
   G_1(x^1,x^2)=1-x^1-x^2, \quad
    G_2(x^1,x^2)=x^1-x^2, \quad
     G_3(x^1,x^2)=x^2.
\]
Here, assumptions (C1), (C2) are fulfilled, the GNEP feasible set $X$ is bounded and satisfies GNEP-Slater. However, (C3) does not hold, since the matrix
\[
   D_x \mathcal{G}_f(x,r) = \left(
  \begin{array}{cc}
      0&  r^1  \\
      \frac{1}{2}r^2 & 0
  \end{array}
  \right)
\]
is indefinite for any parameters $r>0$. Theorem \ref{thm:ex-unique} is thus not applicable. Let us try to apply Theorem \ref{thm:unique} instead, in order to address the issue of local uniqueness.

It is straightforward to see that the Nash equilibria of the given GNEP are of the form 
\[
\bar x(t)=(1-t,t), t \in \left[0,\frac{1}{2}\right].
\] 
The Nash equlibrium $\bar x\left(0\right)=\left(1,0\right)$ is normalized in association with the parameters $r=(r^1,r^2)$, fulfilling $\frac{r^1}{r^2}> \frac{1}{2}$.
The corresponding multipliers for the active constraints from $\mathcal{J}_0(\bar x\left(0\right))=\left\{1,3\right\}$ are
\[
  \bar \Lambda^1= r^1, \quad
  \bar \Lambda^3=r^1-\frac{1}{2}r^2.
\]
Obviously, GNEP-LICQ is fulfilled at $\bar x\left(0\right)$, i.e. also ND1. Due to the choice of the parameters $r$, ND2 holds at $\bar x\left(0\right)$ as well. Since $\mathcal{T}_{\bar x\left(0\right)} X_0(\bar x\left(0\right))=\{0\}$, ND3 is trivially satisfied. The normalized Nash equilibrium $\bar x\left(0\right)$ is thus nondegenerate.

The Nash equlibrium $\bar x\left(\frac{1}{2}\right)=\left(\frac{1}{2},\frac{1}{2}\right)$ is normalized in association with the parameters $r=(r^1,r^2)$, fulfilling $\frac{r^1}{r^2}< \frac{3}{2}$.
The corresponding multipliers for the active constraints from $\mathcal{J}_0(\bar x\left(\frac{1}{2}\right))=\left\{1,2\right\}$ are

\[
  \bar \Lambda^1=\frac{1}{2} r^1+\frac{3}{4} r^2, \quad
  \bar \Lambda^2=-\frac{1}{2}r^1+\frac{3}{4}r^2.
\]
Obviously, GNEP-LICQ is fulfilled at $\bar x\left(\frac{1}{2}\right)$, i.e. also ND1. Due to the choice of the parameters $r$, ND2 holds at $\bar x\left(\frac{1}{2}\right)$ as well. Since $\mathcal{T}_{\bar x\left(\frac{1}{2}\right)} X_0(\bar x\left(\frac{1}{2}\right))=\{0\}$, ND3 is trivially satisfied. The normalized Nash equilibrium $\bar x\left(\frac{1}{2}\right)$ is thus nondegenerate.

The Nash equlibrium $\bar x\left(t\right)=\left(1-t,t\right)$, $t \in \left(0,\frac{1}{2}\right)$, is normalized in association with the parameters $r=(r^1,r^2)$, fulfilling $\frac{r^1}{r^2}= \frac{1+t}{2(1-t)}$.
The corresponding multiplier for the active constraint from $\mathcal{J}_0(\bar x\left(t\right))=\left\{1\right\}$ is
\[
  \bar \Lambda^1= (1-t)r^1= \frac{1+t}{2}r^2.
\]
Obviously, GNEP-LICQ is fulfilled at $\bar x\left(t\right)$, i.e. also ND1. Due to $t \in \left(0,\frac{1}{2}\right)$, ND2 holds at $\bar x\left(0\right)$ as well.
Finally, we have:
\[
  \mathcal{T}_{\bar x} X_0(\bar x) = \left\{(\xi^1,\xi^2)^T \,\left|\, -\xi^1-\xi^2=0\right.\right\}
\]
and
\[
  D_x \mathcal{G}_L(\bar x(t), \bar \Lambda, r) = \left(
  \begin{array}{cc}
      0&  r^1  \\
      \frac{1}{2}r^2 & 0
  \end{array}
  \right).
\] 
Hence, ND3 is also satisfied, since it holds:
\[
 \begin{array}{rcl}
      \displaystyle
  D_x \mathcal{G}_L(\bar x(t), \bar \Lambda, r)\restriction_{\mathcal{T}_{\bar x} X_0(\bar x)} &=& \displaystyle \left(
  \begin{array}{cc}
      1 & -1 
  \end{array}
  \right)  \left(
  \begin{array}{cc}
      0&  r^1  \\
      \frac{1}{2}r^2 & 0
  \end{array}
  \right)
  \left(
  \begin{array}{c}
      1 \\ -1 
  \end{array}
  \right) \\ \\
  &=& \displaystyle - r^1 - 2\varepsilon r^2 \not =0.
  \end{array}
\] 
The normalized Nash equilibrium $\bar x\left(t\right)$, $t \in \left(0,\frac{1}{2}\right)$, is thus nondegenerate. Moreover, we deduce by recalling that $t \in \left(0,\frac{1}{2}\right)$:
\[  
 \frac{1}{2} < \frac{r^1}{r^2} < \frac{3}{2}.
\]
Vice versa, if the parameters $r=(r^1,r^2)$ fulfill these inequalities, there is the unique normalized Nash equilibrium associated with $r$ of the form
\[
   \bar x\left(t\right)=\left(1-t,t\right) \mbox{ with } t=\frac{\frac{r^1}{r^2} - \frac{1}{2}}{\frac{r^1}{r^2} + \frac{1}{2}}.
\]

Overall, we distinguish the following cases with respect to the (local) uniqueness of normalized Nash equilibria:
\begin{center}
    \begin{tabular}{|c|c|}
    \hline
        Parameters $r=\left(r^1,r^2\right)$ & Normalized Nash equilibria associated with $r$\\ \hline
        $0 < \frac{r^1}{r^2} \leq \frac{1}{2}$ & $\bar x\left(\frac{1}{2}\right)$\\ \hline
        $\frac{1}{2} < \frac{r^1}{r^2} < \frac{3}{2}$ & $\bar x\left(0\right)$, $\bar x\left(\frac{1}{2}\right)$, $\bar x\left(t\right)=\left(1-t,t\right) \mbox{ with } t=\frac{\frac{r^1}{r^2} - \frac{1}{2}}{\frac{r^1}{r^2} + \frac{1}{2}}$  \\ \hline
        $\frac{r^1}{r^2} \geq \frac{3}{2}$ & $\bar x\left(0\right)$\\ \hline
    \end{tabular}
\end{center}
Note that in every of the three cases Theorem \ref{thm:unique} is applicable, since all normalized Nash equilibria are nondegenerate. This is confirmed by the observation that the normalized Nash equlibria associated with the same parameters are at least locally unique. We emphasize that the phenomenon presented here is stable with respect to arbitrarily small perturbations of the defining functions. \qed 
\end{example}

\section*{Acknowledgements}
The author would like to thank Bernd Kummer for a valuable discussion on the notion of normalized Nash equilibrium.


\begin{thebibliography}{99}

\bibitem{dorsch:2013} D. Dorsch, H. Th. Jongen, V. Shikhman, {\it On structure and computation of generalized Nash equilibria},
SIAM Journal on Optimization {\bf 23} (2013), 452--474.

\bibitem{dorsch:2013a} D. Dorsch, H. Th. Jongen, V. Shikhman, {\it On intrinsic complexity of Nash equilibrium problems and bilevel optimizatio},
Journal of Optimization Theory and Applications {\bf 159} (2013), 606--634.

\bibitem{faccinei:2007} F. Facchinei, A. Fischer, V. Picciallia, {\it On generalized Nash games and variational inequalities},
Operations Research Letters {\bf 35} (2007), 159--164.

\bibitem{faccinei:2010} F. Facchinei, C. Kanzow, {\it Generalized Nash equilibrium problems},
Annals of Operations Research {\bf 175} (2010), 177--211.

\bibitem{fukushima:2011} M. Fukushima, {\it Restricted generalized Nash equilibria
and controlled penalty algorithm},
Computational Management Science {\bf 8} (2011), 201--218.

\bibitem{guenzel:2008} H. G\"unzel, {\it The structured jet transversality theorem},
Optimization {\bf 57} (2008), 159--164.

\bibitem{heusinger:2012} A. von Heusinger, C. Kanzow, M. Fukushima, {\it Newton’s method for computing a normalized equilibrium in the generalized Nash game
through fixed point formulation},
Mathematical Programming {\bf 132} (2012), 99--123.

\bibitem{hirsch:1976} M. W. Hirsch, 
	{\it Differential Topology}, Springer, Berlin-Heidelberg-New York, 1976.

\bibitem{jongen:2000} H. Th. Jongen, P. Jonker, F. Twilt, 
	{\it Nonlinear Optimization in Finite Dimensions}, Kluwer Academic Publishers, Dordrecht, 2000.

\bibitem{jongen:2004} H. Th. Jongen, K. Meer, E. Triesch, 
	{\it Optimization Theory}, Kluwer Academic Publishers, Dordrecht, 2004.


\bibitem{kakutani:1941} S. Kakutani, {\it A generalization of Brouwer's fixed point theorem},
Duke Mathematics Journal {\bf 8} (1941), 457--459.

\bibitem{nabetani:2011} K. Nabetani, P. Tseng, M. Fukushima, {\it Parametrized variational inequality approaches to generalized Nash equilibrium problems with shared constraints},
Computational Optimization and Applications {\bf 48} (2011), 423--452.

\bibitem{rosen:1965} J. B. Rosen, {\it Existence and uniqueness of equilibrium points for concave N-person games},
Econometrica {\bf 33} (1965), 520--534.

\bibitem{stein:2018} O. Stein, 
	{\it Grundz\"uge der
Nichtlinearen Optimierung}, Springer Spektrum, Berlin, 2018.
	
\end{thebibliography}
\end{document}